\documentclass[12pt]{article}
\usepackage{amssymb}
\usepackage{amsmath}
\usepackage{amscd}
\usepackage[all]{xy}
\usepackage[french]{babel}

\usepackage{color}
\title{Quelques cas d'annulation du troisi\`eme groupe de cohomologie non ramifi\'ee}
\author{Jean-Louis Colliot-Th\'el\`ene }
\date{6 octobre 2011}

\def\beq{\begin{equation} \label}

\newtheorem{theo}{Th\'eor\`eme}[section]
\newtheorem{prop}[theo]{Proposition}
\newtheorem{lem}[theo]{Lemme}
\newtheorem{cor}[theo]{Corollaire}
\newtheorem{defi}[theo]{D\'efinition}

\newcommand{\bthe}{\begin{theo}}
\newcommand{\ble}{\begin{lem}}
\newcommand{\bpr}{\begin{prop}}
\newcommand{\bco}{\begin{cor}}
\newcommand{\bde}{\begin{defi}}
\newcommand{\ethe}{\end{theo}}
\newcommand{\ele}{\end{lem}}
\newcommand{\epr}{\end{prop}}
\newcommand{\eco}{\end{cor}}
\newcommand{\ede}{\end{defi}}

\newcommand{\oi}{\hskip1mm {\buildrel \simeq \over \rightarrow} \hskip1mm}

\def \to {{\rightarrow}}

\def \Q {{\mathbb Q}}
\def \N {{\mathbb N}}
\def \Z {{\mathbb Z}}
\def \F {{\mathbb F}}
\def \C {{\mathbb C}}

\def \P {{\mathbf P}}

\begin{document}

\maketitle

{\bf R\'esum\'e} On \'etablit la nullit\'e du troisi\`eme groupe de cohomologie non ramifi\'ee
pour certaines vari\'et\'es munies d'un pinceau   de quadriques ou d'intersections compl\`etes lisses de deux quadriques.  
Sur les complexes, ceci permet d'\'etablir la conjecture de Hodge enti\`ere en
degr\'e 4 pour de telles vari\'et\'es.
\noindent  

\medskip

{\bf Abstract}  The third unramified cohomology group is shown to vanish on certain varieties
equipped with a pencil of quadrics or  of  smooth complete intersections of two quadrics. Over the complex field,
this shows that the integral Hodge conjecture in degree 4  holds for such  varieties.

MSC-class :  14C35; 14E08,   14F99

\section{Notations et rappels}

Cet article apporte un compl\'ement  \`a des textes  r\'ecents de C.~Voisin et l'auteur \cite{ctv}, de B.~Kahn et l'auteur  \cite{CTK} et  de C.~Voisin \cite{V}.
Les notations sont celles de ces articles.
La cohomologie des corps ici utilis\'ee est la cohomologie galoisienne,
dont nous utilisons librement les propri\'et\'es \cite{serre}.
Pour $F$ un corps et $M$ un module galoisien discret sur le groupe de Galois absolu de $F$,
et $i \geq 0$ un entier, on note $H^{i}(F,M)$ le $i$-\`eme groupe de cohomologie galoisienne
\`a valeurs dans $M$.
Pour les propri\'et\'es de  la cohomologie non ramifi\'ee des vari\'et\'es et de leurs corps
de fonctions, on consultera  \cite{CTBarbara}.
Pour $X$ une vari\'et\'e int\`egre sur un corps $k$, on note $k(X)$ son corps des fonctions.
Pour $i \in \N$ et $j \in \Z$, et $l$ un nombre premier distinct de la caract\'eristique de $k$, 
on note $\Q_{l}/\Z_{l}(j)$ la limite directe  sur $n$ des groupes de racines de l'unit\'e tordus 
$\mu_{l^n}^{\otimes j}$, et on note
$$H^{i}_{nr}(k(X)/k,\Q_{l}/\Z_{l}(j)) \subset H^{i}(k(X),\Q_{l}/\Z_{l}(j))$$
le sous-groupe consistant des \'el\'ements $\xi \in H^{i}(k(X),\Q_{l}/\Z_{l}(j))$  dont tous les r\'esidus $\delta_{A}(\xi)$ sont nuls,
o\`u $A$ parcourt l'ensemble des anneaux de valuation discr\`ete de rang 1 contenant $k$ et
de corps des fractions $k(X)$.

\begin{lem}\label{lemmetrivial}
Soit $F$ un corps de caract\'eristique diff\'erente de 2 et soit $Z$ une $F$-quadrique lisse
de dimension au moins 1.
Pour tout $l$ premier impair distinct de la caract\'eristique, pour tout entier $i \geq 0$,
et pour tout $j \in \Z$, l'application de restriction naturelle
$H^{i}(F,\Q_{l}/\Z_{l}(j)) \to H^{i}_{nr}(F(Z),\Q_{l}/\Z_{l}(j))$ est un isomorphisme.
\end{lem}
{\it D\'emonstration}. Il existe une extension quadratique s\'eparable  $L/F$ sur laquelle
$Z$ acquiert un $L$-point, donc est $L$-birationnelle \`a un espace projectif.
Ceci implique $H^{i}(L,\Q_{l}/\Z_{l}(j)) \simeq  H^{i}_{nr}(L(Z),\Q_{l}/\Z_{l}(j))$,
et le r\'esultat annonc\'e s'obtient par un argument de trace. QED

\begin{prop}\label{ctsasd}(\cite[Thm. 3.2, p.~60]{CSS})
Soit $F$ un corps de caract\'eristique diff\'erente de 2. Soit $X \subset \P^n_{F}$
une intersection compl\`ete lisse de deux quadriques.  Supposons $n\geq 4$.
Si $X$ poss\`ede un point rationnel sur $F$, alors $X$ est $F$-birationnelle \`a
une $F$-vari\'et\'e g\'eom\'etriquement int\`egre $Z$ munie d'un morphisme
$Z \to \P^1_{F}$ dont la fibre g\'en\'erique est une quadrique lisse dans $\P^{n-2}_{F(\P^1)}$.
\end{prop}

\begin{theo}\label{karosu} 
Soit $F$ un corps de   caract\'eristique diff\'erente de 2. Soit $Z$ une $F$-quadrique lisse
de dimension $d\geq 1$ qui n'est pas une quadrique d'Albert anisotrope.
Pour tout $l$ premier diff\'erent de la caract\'eristique,
l'application  de restriction naturelle sur la cohomologie galoisienne
$$H^3(F,\Q_{l}/\Z_{l}(2)) \to H^3(F(Z),\Q_{l}/\Z_{l}(2))$$
induit une surjection
$$H^3(F,\Q_{l}/\Z_{l}(2)) \to H^3_{nr}(F(Z),\Q_{l}/\Z_{l}(2))$$
\end{theo}

Une quadrique d'Albert est une quadrique de dimension 4
d\'efinie par une forme quadratique diagonale $<a,b,-ab,-c,-d,cd>$, avec $a,b,c,d \in F^{\times}$.

Pour $l\neq 2$, cela r\'esulte du lemme \ref{lemmetrivial}.
Pour $l=2$,  cas o\`u
il y a  aussi des variantes plus d\'elicates avec les coefficients $\Z/2$,
plusieurs auteurs (Suslin, Merkur'ev, Peyre, Sujatha, Kahn, Rost) 
ont contribu\'e au th\'eor\`eme ci-dessus.
L'\'enonc\'e g\'en\'eral fait l'objet de \cite[\S 10]{KRS}, article auquel on
se r\'ef\'erera pour le d\'etail des diff\'erentes contributions au sujet.

{\it Remarque}
Comme il est expliqu\'e dans \cite[Prop. 3]{PirutkaCRAS}, gr\^ace \`a des travaux de B.~Kahn, 
on dispose d'un \'enonc\'e analogue au th\'eor\`eme \ref{karosu} pour les vari\'et\'es
de Severi-Brauer d'indice premier. Ceci permet d'\'etablir pour les fibrations en de telles
vari\'et\'es les analogues de tous
les \'enonc\'es donn\'es ici pour les fibrations en coniques. 
 
\section{Le  th\'eor\`eme g\'en\'eral}

\begin{theo}\label{principal}
Soit $k$ un corps de   caract\'eristique diff\'erente de 2 
et soit $l$ un nombre premier distinct de la caract\'eristique de $k$.
Soit $Y$ une $k$-vari\'et\'e g\'eom\'etriquement int\`egre.
Soit $X$ une $k$-vari\'et\'e connexe projective et lisse de dimension  $d$
qui est $k$-birationnelle \`a une $k$-vari\'et\'e g\'eom\'etriquement int\`egre $Z$
munie d'un $k$-morphisme dominant $Z \to Y$ de fibre g\'en\'erique
une quadrique lisse $Q$ sur le corps $k(Y)$. Supposons $dim(Q)= d-dim(Y) \geq 1$.

Dans chacun des cas suivants :

\smallskip

(a)  $k$ est s\'eparablement clos  et $dim(Y) \leq 2$,

(b) $k$ est un corps de $l$-dimension cohomologique 1, 
 $Q$ n'est pas une quadrique d'Albert anisotrope, et   $dim(Y) \leq 1$,
 
(c) $k$ est un corps $C_{1}$  et   $dim(Y) \leq 1$,

\smallskip

\noindent on a $H^3_{nr}(k(X)/k,\Q_{l}/\Z_{l}(2))=0$,

\end{theo}

{\it D\'emonstration} 
 Supposons d'abord $dim(Q) \geq 3$.

Si $k$ est un corps $C_{1}$ et $dim(Y) \leq 1$, alors $k(Y)$ est un corps $C_{2}$.
La  quadrique $Q$, de dimension au moins 3, poss\`ede donc  un point rationnel sur $k(Y)$,
donc, comme elle est lisse, est $k(Y)$-birationnelle \`a un espace projectif sur $k(Y)$.

Supposons $k$  s\'eparablement clos et  $dim(Y) \leq 2$.
 Soit $\overline{k}$ une cl\^oture alg\'ebrique de $k$.
Le corps $\overline{k}(Y)$ est un corps $C_{2}$, donc  la
quadrique $Q$, de dimension au moins $3$, sur $\overline{k}(Y)$ poss\`ede un point rationnel.
Comme l'extension $\overline{k}(Y)/k(Y)$ est de pro-degr\'e impair, 
un th\'eor\`eme  bien connu de T. A. Springer implique que   $Q$ 
 poss\`ede un point rationnel sur $k(Y)$, donc comme  elle est lisse, est $k(Y)$-birationnelle
\`a un espace projectif. 

On a donc dans ces deux cas
$$ H^3(k(Y), \Q_{l}/\Z_{l}(2) )\oi H^3_{nr}(k(Y)(Q)/k(Y),\Q_{l}/\Z_{l}(2)).$$
 
\medskip

Supposons $dim(Q) \leq 2$. La quadrique $Q$ n'est donc pas
une quadrique d'Albert.
D'apr\`es le th\'eor\`eme \ref{karosu},  l'application de restriction
$$H^3(k(Y),\Q_{l}/\Z_{l}(2)) \to H^3_{nr}(k(Y)(Q)/k(Y),\Q_{l}/\Z_{l}(2))$$
est surjective.

\medskip

   Sous les hypoth\`eses du th\'eor\`eme,
la $l$-dimension cohomologique de $k(Y)$ est au plus 2. On a donc
$H^3(k(Y),\Q_{l}/\Z_{l}(2))=0$. D'apr\`es ce qui pr\'ec\`ede, on a donc
   $H^3_{nr}(k(Y)(Q)/k(Y),\Q_{l}/\Z_{l}(2))=0$.  On a $k(X)=k(Y)(Q)$ et  l'inclusion
 $$H^3_{nr}(k(X)/k,\Q_{l}/\Z_{l}(2)) \subset H^3_{nr}(k(Y)(Q)/k(Y),\Q_{l}/\Z_{l}(2)).$$
 On a donc bien  $H^3_{nr}(k(X)/k,\Q_{l}/\Z_{l}(2))=0$.
 QED

\section{Sur les complexes}

\begin{cor}\label{corcomplexe} 
Soit $f : X \to Y$ un morphisme dominant de vari\'et\'es connexes projectives 
et lisses sur le corps $\C$ des complexes.  Dans chacun des cas suivants :

\smallskip

(a) $dim(Y) \leq 2$ et la fibre g\'en\'erique de $ f$ est une quadrique lisse de dimension
au moins 1,

(b) $Y=\Gamma$ est une courbe et la fibre g\'en\'erique de $f$ est une intersection compl\`ete  lisse de deux
quadriques de dimension au moins 2 sur le corps $\C(\Gamma)$,  

\smallskip
\noindent on a $H^3_{nr}(\C(X)/\C,\Q/\Z(2))=0$, et  la conjecture de Hodge enti\`ere vaut pour les
classes de Hodge enti\`eres de degr\'e 4 sur $X$ : toute telle classe dans $H^4_{\rm Betti}(X,\Z)$ est l'image d'un cycle  alg\'ebrique de codimension 2 sur $X$.
\end{cor}

{\it D\'emonstration}
Soit $Z^4(X)$ le quotient du groupe des classes de Hodge dans $H^4_{\rm Betti}(X,\Z)$
par le sous-groupe des classes de cycles alg\'ebriques.

Dans chacun des deux cas consid\'er\'es, le groupe de Chow 
des z\'ero-cycles de degr\'e z\'ero sur $X$  est support\'e sur une surface.
Le th\'eor\`eme 1.1 de \cite{ctv} donne   $H^3_{nr}(\C(X)/\C,\Q/\Z(2))\simeq Z^4(X)$.

Dans le cas (a), le th\'eor\`eme  \ref{principal}(a)  ci-dessus donne  $$H^3_{nr}(\C(X)/\C,\Q/\Z(2))=0.$$

Dans le cas (b), comme $\C(\Gamma)$ est un corps $C_{1}$,  la fibre g\'en\'erique
de $f$ poss\`ede un point rationnel sur $\C(\Gamma)$. La  proposition \ref{ctsasd}
montre  alors que cette fibre g\'en\'erique est $\C(\Gamma)$-birationnelle
\`a une fibration en quadriques $Z \to \P^1_{\C(\Gamma)}$
 de dimension relative au moins 1 au-dessus de $\P^1_{\C(\Gamma)}$, \`a fibre g\'en\'erique lisse.
 La $\C$-vari\'et\'e $X$ est donc birationnelle \`a l'espace total d'une fibration $Z \to S= \P^1\times \Gamma$
 dont la fibre g\'en\'erique est une quadrique lisse $Q/ \C(S)$ de dimension au moins 1.
 On a l'inclusion
 $$H^3_{nr}(\C(X)/\C,\Q/\Z(2)) \subset  H^3_{nr}(\C(S)(Q) /\C(S),\Q/\Z(2))$$
 Les th\'eor\`emes  \ref{karosu} et \ref{principal}  montrent que ce dernier groupe
  est un quotient de 
 $H^3(\C(S),\Q/\Z(2))$, qui est nul car $\C(S)$ est de dimension cohomologique 2. QED
 \medskip

{\it Commentaires} 
  
Le cas (a) est mis ici pour m\'emoire, il a d\'ej\`a \'et\'e \'etabli dans \cite[Cor. 8.2]{ctv}.

Dans le cas (b), pour une famille  de dimension relative $d=2$, l'\'enonc\'e
est un cas particulier de \cite[Thm. 8.14]{ctv}. Pour une famille de dimension relative
$d=3$,  l'\'enonc\'e obtenu g\'en\'eralise le Cor. 1.7  de \cite{V}, \'etabli par C. Voisin par des m\'ethodes
g\'eom\'etriques, au prix d'hypoth\`eses restrictives sur les fibres singuli\`eres
de $X \to \Gamma$. Ici on ne fait aucune hypoth\`ese sur ces fibres. Le cas de la dimension relative
$d \geq 4$ est facile : une quadrique lisse de dimension au moins 3 sur $\C(\Gamma)$
poss\`ede un point rationnel sur $\C(\Gamma)$, donc est $\C(\Gamma)$-birationnelle
\`a un espace projectif sur $\C(\Gamma)$. La vari\'et\'e $X$ est donc
 $\C$-birationnel \`a un produit $\P^{d} \times \Gamma$. 
Pour une fibration de $\C(\Gamma)$-vari\'et\'es $Z \to \P^1_{\C(\Gamma)}$
de fibre g\'en\'erique une quadrique lisse de dimension 1 ou 2, la $\C(\Gamma)$-vari\'et\'e
$Z$ n'est pas n\'ecessairement $\C(\Gamma)$-birationnelle
\`a un espace projectif sur $\C(\Gamma)$, comme on voit en calculant le
groupe de Brauer non ramifi\'e de $Z$.
 Mais c'est  semble-t-il une question ouverte si une intersection compl\`ete  lisse de deux quadriques
 dans  $\P^5_{\C(\Gamma)}$ est $\C(\Gamma)$-birationnelle \`a un espace projectif sur
 $\C(\Gamma)$.

Soit maintenant $X$ une vari\'et\'e connexe, projective et lisse munie d'une fibration $X \to \Gamma$
sur une courbe $\Gamma$, dont la fibre g\'en\'erique est une hypersuface cubique dans $\P^{d+1}_{\C(\Gamma)}$.
Pour $d=2$, on peut, soit par la g\'eom\'etrie complexe (Voisin, cf. \cite[Thm. 6.1]{ctv} ou \cite[Thm. 1.3]{V}),  soit par la $K$-th\'eorie  alg\'ebrique \cite[Thm. 8.14]{ctv},  montrer   que la conclusion du corollaire ci-dessus vaut encore. 

Pour $d=3$, sous des hypoth\`eses restrictives sur les fibres singuli\`eres
de $X \to \Gamma$,  ceci vaut encore pour $X$. C'est l\`a un r\'esultat r\'ecent
de C.~Voisin \cite[Thm. 2.11]{V}.
Ce r\'esultat, obtenu par des m\'ethodes g\'eom\'etriques, semble  hors d'atteinte des m\'ethodes  
de $K$-th\'eorie alg\'ebrique.

 \section{Sur les corps finis}

\begin{cor}\label{corfini} 
Soit $\F$ un corps fini de caract\'eristique diff\'erente de 2
et soit $l$ premier distinct de la caract\'eristique de $\F$.
 Soit $X$ une $\F$-vari\'et\'e projective, lisse, g\'eom\'etriquement int\`egre.
 Dans chacun des cas suivants :

(a) il existe un $\F$-morphisme dominant $f: X \to C$  de  $X$ vers  une $\F$-courbe 
$C$ projective, lisse et g\'eom\'etriquement int\`egre,
et la fibre g\'en\'erique de $f$  est une
quadrique lisse de dimension au moins 1 sur $\F(C)$,

(b) la $\F$-vari\'et\'e $X$ est une intersection  compl\`ete lisse de deux
quadriques dans $\P^n_{\F}$ et $n \geq 4$,

\noindent  on a
$H^3_{nr}(\F(X)/\F,\Q_{l}/\Z_{l}(2))=0$.
 
\end{cor}

{\it D\'emonstration}
L'\'enonc\'e (a) est une application imm\'ediate du th\'eor\`eme \ref{principal}(c), puisqu'un corps fini
est un corps $C_{1}$. Comme $\F$ est $C_{1}$, toute vari\'et\'e $X$ comme en (b) admet un
point $\F$-rationnel. La proposition \ref{ctsasd} montre que $X$ est $\F$-birationnelle \`a 
une $\F$-vari\'et\'e g\'eom\'etriquement int\`egre $Z$ munie d'un $\F$-morphisme dominant
$Z \to \P^1_{\F}$ de fibre g\'en\'erique une quadrique lisse de dimension au moins 1.
On a l'inclusion $$H^3_{nr}(\F(X)/\F,\Q_{l}/\Z_{l}(2)) \subset H^3_{nr}(\F(Z)/\F(\P^1), \Q_{l}/\Z_{l}(2)).$$
Le  th\'eor\`eme \ref{karosu} montre que l'application de restriction 
$$H^3(\F(\P^1) ,  \Q_{l}/\Z_{l}(2)) \to
    H^3_{nr}(\F(Z)/\F(\P^1), \Q_{l}/\Z_{l}(2))$$ est surjective. Mais $H^3(\F(\P^1) ,  \Q_{l}/\Z_{l}(2)) =0$,
    car le corps 
$\F(\P^1)$ est de $l$-dimension cohomologique 2.
    QED

\medskip

{\it Commentaires}

Soit $X$ une vari\'et\'e projective, lisse, g\'eom\'etriquement int\`egre sur un corps $\F$.
Soit $l$ un premier distinct de la caract\'eristique de $\F$.

C'est un th\'eor\`eme de th\'eorie
du corps de classes sup\'erieur (cf. \cite[Prop. 3.1]{CTK} 
que pour $X$ de dimension 2, on a $H^3_{nr}(\F(X),\Q_{l}/\Z_{l}(2))=0$.

En dimension au moins 3,  comme il est expliqu\'e dans  \cite{CTK},
la situation est la suivante.
On conjecture que  $H^3_{nr}(\F(X),\Q_{l}/\Z_{l}(2))$ est un groupe fini.
On ne conna\^{\i}t pas une seule vari\'et\'e $X/\F$ de dimension au plus 4
avec $H^3_{nr}(\F(X),\Q_{l}/\Z_{l}(2))Ê\neq 0$. 
On conna\^it des vari\'et\'es de dimension 5 sur un corps fini, fibr\'ees en quadriques
de dimension 3 au-dessus d'une surface, donc 
g\'eom\'etriquement rationnelles, avec $H^3_{nr}(\F(X),\Q_{l}/\Z_{l}(2))Ê\neq 0$
(Pirutka \cite{Pirutka}).
On conjecture que l'on a $H^3_{nr}(\F(X),\Q_{l}/\Z_{l}(2))Ê= 0$
 pour les vari\'et\'es de dimension 3 qui sont g\'eom\'etriquement unir\'egl\'ees,
 par exemple les vari\'et\'es fibr\'ees en surfaces cubiques au-dessus d'une courbe.
 
Supposons la caract\'eristique de $\F$ impaire.
Un th\'eor\`eme d\'elicat de Parimala et Suresh (cf. \cite[Thm. 4.4]{CTK}) \'etablit
$H^3_{nr}(\F(X),\Q_{l}/\Z_{l}(2))Ê= 0$  
pour les vari\'et\'es de dimension 3 fibr\'ees en coniques au-dessus d'une surface.
Le corollaire \ref{corfini} montre que l'on a  $H^3_{nr}(\F(X),\Q_{l}/\Z_{l}(2))Ê= 0$
pour les vari\'et\'es   fibr\'ees en quadriques  de dimension au moins 1  au-dessus d'une courbe,
et aussi pour les intersections compl\`etes lisses de deux quadriques dans $\P^n_{\F}$ pour $n \geq 5$.
Ces derni\`eres sont  des vari\'et\'es  g\'eom\'e\-tri\-quement rationnelles. Elle sont en fait
$\F$-birationnelles \`a un espace projectif sur leur corps de base $\F$, ce qui donne une d\'emonstration directe du Corollaire \ref{corfini} (b). 
C'est  facile \`a voir pour
$n \geq 6$ (\cite[Thm. 3.4]{CSS}), car dans ce cas, $X$ est birationnelle \`a une famille
de quadriques de dimension au moins 3 sur $\P^1_{\F}$, et $\F(\P^1)$ est un corps $C_{2}$.
Le cas $n=5$ est une cons\'equence d'un fait relativement peu connu : 
une  
intersection compl\`ete lisse $X$ de deux quadriques dans $\P^5_{\F}$ 
contient une droite, d\'efinie sur $\F$, de $\P^n_{\F}$, 
ce qui implique  (cf. \cite[Prop. 2.2]{CSS}) que  $X$ est 
$\F$-birationnelle \`a un espace projectif.  
La m\'ethode donn\'ee ici pour \'etablir (b) est n\'eanmoins int\'eressante,
elle s'applique encore \`a certaines  intersections compl\`etes singuli\`eres
de deux quadriques.

\vskip2cm

\noindent
CNRS, UMR 8628, Math\'ematiques, B\^atiment 425, 

\noindent
Universit\'e Paris-Sud, F-91405 Orsay, France
\medskip

\noindent jlct@math.u-psud.fr


\begin{thebibliography}{99}

\bibitem{CTBarbara}    J.-L. Colliot-Th\'{e}l\`{e}ne,
 Birational invariants, purity and the Gersten conjecture,
in  K-Theory and Algebraic Geometry: Connections with Quadratic Forms and Division Algebras, AMS Summer Research Institute, Santa Barbara 1992, ed. W. Jacob and A. Rosenberg, Proceedings of Symposia in Pure Mathematics {\bf  58}, Part I (1995) p. 1--64.

\bibitem{CSS}  J.-L. Colliot-Th\'{e}l\`{e}ne, J.-J. Sansuc et  Sir Peter Swinnerton-Dyer, 
Intersection of two quadrics and Ch\^atelet surfaces, I, Crelle {\bf 373} (1987) 37-107.

\bibitem{CTK} J.-L. Colliot-Th\'{e}l\`{e}ne et B. Kahn, 
Cycles de codimension 2 et $H^3$ non ramifi\'e pour les vari\'et\'es sur les corps finis, 
{\tt{ arXiv:1104.3350v1 [math.AG]}}.

\bibitem{ctv}  J.-L. Colliot-Th\'{e}l\`{e}ne et C. Voisin, Cohomologie non ramifi\'ee et
conjecture de Hodge enti\`ere,  
{\tt arXiv:1005.2778v2 [math.AG]}, \`a para\^{\i}tre dans Duke Mathematical Journal.

\bibitem{KRS} B. Kahn, M. Rost et R. Sujatha, Unramified cohomology of quadrics. I, Amer. J. Math. {\bf  120} (1998) 841--891.

\bibitem{Pirutka} A. Pirutka,
Sur le groupe de Chow de codimension deux des vari\'et\'es  sur les corps finis  
{{\tt : arXiv :math/1004.1897v2}}, \`a para\^{\i}tre dans Algebra  and Number Theory. 

\bibitem{PirutkaCRAS} A. Pirutka, Cohomologie non ramifi\'ee en degr\'e trois d'une vari\'et\'e
de Severi-Brauer, 
C. R. Acad. Sci. Paris, S\'er. I  {\bf 349} (2011) 369Ð373 



\bibitem{serre} J-P. Serre, Cohomologie galoisienne,  5\`eme \'ed.,  r\'evis\'ee et compl\'et\'ee,
 Springer LNM {\bf 5} (1994).

 
\bibitem{V} C. Voisin,  Abel-Jacobi maps, integral Hodge classes and decomposition of the diagonal,  {\tt{arXiv:1005.5621v2 [math.AG]}}.


\end{thebibliography}
\end{document}